# ESTIMATING THE NUMBER OF CLASSES

By Chang Xuan Mao and Bruce G. Lindsay

*University of California, Riverside and Pennsylvania State University*

Estimating the unknown number of classes in a population has numerous important applications. In a Poisson mixture model, the problem is reduced to estimating the odds that a class is undetected in a sample. The discontinuity of the odds prevents the existence of locally unbiased and informative estimators and restricts confidence intervals to be one-sided. Confidence intervals for the number of classes are also necessarily one-sided. A sequence of lower bounds to the odds is developed and used to define pseudo maximum likelihood estimators for the number of classes.

**1. Introduction.** The species problem has a wide variety of applications [3]. The term "species" has been endowed with many meanings such as taxa, words known by an author and expressed genes in a tissue. Consider a population of infinitely many individuals belonging to $c$ distinct classes labeled by $i = 1, 2, \ldots, c$. In a sample of $S$ individuals, $Y_i$ individuals belong to the $i$th class. The $i$th class is not detected when $Y_i = 0$. Estimating the number of classes $c$ from those $Y_i > 0$ is a well-known difficult problem. For example, I. J. Good pointed out that "*I don't believe it is usually possible to estimate the number of species, but only an appropriate lower bound to that number. This is because there is nearly always a good chance that there are a very large number of extremely rare species*" [3].

In the literature, $Y_i$ is usually treated as a Poisson random variable with mean $\lambda_i$ and the $\lambda_i$ are assumed to follow a mixing distribution $P$ over $(0, \infty)$. The $Y_i$ arise as a sample from $g_P(y) = \int e^{-\lambda} \lambda^y/(y!) \, dP(\lambda)$. There are $n = \sum_{i=1}^{c} I(Y_i > 0)$ detected classes in the sample. Because $n \sim$ binomial$(c, 1 - g_P(0))$, the maximum likelihood estimator (MLE) of $c$ given $P$ is the integer part of $\hat{c}(\theta) = n(1 + \theta)$, where $\theta = g_P(0)/\{1 - g_P(0)\}$ is the odds that a single class is undetected in the sample. When an estimator $\hat{\theta}$ for $\theta$ is substituted into $\hat{c}(\theta)$, we obtain a pseudo MLE for $c$ [11]. The problem of estimating $c$ is thereby reduced to that of estimating $\theta$.

The idea of the nonexistence of inferential procedures is not unfamiliar to statisticians (e.g., [1, 10, 13, 17]). To make Good's point concrete, we will show that the odds $\theta$ is discontinuous. There are several consequences: no locally unbiased and informative estimator for $\theta$, no genuine two-sided confidence intervals and arbitrarily bad informativity when reducing bias to zero. However, because $\theta$ is lower semicontinuous, there exist nonparametric lower confidence limits. A sequence of closed-form lower bounds to $\theta$ is developed. Similar results concerning inference on the number of classes $c$ hold. The upper confidence limits for $c$ are often infinite. The estimators for the lower bounds to $\theta$ yield estimators for lower bounds to $c$.

This article is organized as follows. The mixture model will be described in Section 2. In Section 3 the discontinuity of $\theta$ and its consequences will be investigated. In Section 4 we will demonstrate the lower semicontinuity of $\theta$, construct lower confidence limits and propose its lower bounds. The problem of estimating $c$ will be considered in Section 5. In Section 6 an epidemiological application and a genomic application will be studied. In Section 7 extensions to related problems will be discussed. All proofs are contained in Section 8.

**2. The mixture model.** Let $n_y = \sum_{i=1}^{c} I(Y_i = y)$. Since the $Y_i$ arise from $g_P(y)$, $(n_0, n_1, \ldots)$ follows a multinomial density. When $n_0$, the number of classes in the population unobserved in the sample, is replaced with $c - n$, this yields

$$p_1(c, P) = \frac{c!}{(c-n)! \prod_{x=1}^{\infty} n_x!} g_P^{c-n}(0) \prod_{x=1}^{\infty} g_P^{n_x}(x).$$

This likelihood can be written as $p_1(c, P) = p_2(c, P) p_3(n, P)$, where $p_2(c, P)$ is the density of $n$ and $p_3(n, P)$ is the conditional density of $(n_1, n_2, \ldots)$ given $n$:

$$p_2(c, P) = \binom{c}{n} g_P^{c-n}(0) \{1 - g_P(0)\}^n,$$

$$p_3(n, P) = \frac{n!}{\prod_{x=1}^{\infty} n_x!} \prod_{x=1}^{\infty} \left\{ \frac{g_P(x)}{1 - g_P(0)} \right\}^{n_x}.$$

The likelihood of $n$ is binomial, as indicated before, and depends on both $c$ and $\theta$. The conditional likelihood has no dependence on $c$, but contains most of the information about $P$. It can be analyzed as follows. Conditioning on $n$, those $Y_i > 0$ follow a zero-truncated mixture $g_P(x)/\{1 - g_P(0)\}$. Rewrite them as $X_1, X_2, \ldots, X_n$. Let $f_\lambda(x) = \lambda^x / \{x!(e^\lambda - 1)\}$ and $f_Q(x) = \int f_\lambda(x) \, dQ(\lambda)$, where

(2.1) $$dQ(\lambda) = (1 - e^{-\lambda}) \, dP(\lambda) \bigg/ \int (1 - e^{-\lambda}) \, dP(\lambda).$$

Because $f_Q(x) = g_P(x)/\{1 - g_P(0)\}$ for $x \geq 1$, the $X_i$ can be treated as a sample from a mixture of zero-truncated Poisson densities. The joint density of the $X_i$ is

$$f_Q^{(n)}(x_1, x_2, \ldots, x_n) = \prod_{i=1}^{n} f_Q(x_i) = \prod_{x=1}^{\infty} f_Q^{n_x}(x).$$

Note that $Q$ has no mass on zero. The nonparametric mixture model refers to $\mathcal{F} = \{f_Q : Q \in \mathcal{Q}\}$, where $\mathcal{Q}$ contains all legitimate mixing distributions.

LEMMA 2.1. *$\mathcal{F}$ is identifiable in the sense that $f_Q = f_G$ yields $Q = G$.*

Finally, note that $n$ plays a dual role as the number of detected classes and the sample size of the $X_i$, and $c$ also plays a dual role as the parameter of interest and the "sample size" of the $Y_i$. Our asymptotic results concerning $\theta$-estimation will be based on $n$ becoming infinite, which implies that $c$ goes to infinity as $c = E(n) \times (1 + \theta)$, a common natural practice in the literature that deals with nonstandard problems with integer parameters. However, our key result for $c$-estimation will be finite-sample in nature, so that no asymptotics are required.

**3. Discontinuity.** We will show that estimating $\theta$ is difficult in several aspects.

We write $\theta = \theta(f_Q)$ because of Lemma 2.1 and the fact that $\theta = \int (e^{\lambda} - 1)^{-1} dQ(\lambda)$. As the mass of $Q$ at zero is nonidentifiable and mass near zero is nearly undetectable, we have the following result.

LEMMA 3.1. *$\theta$ is Hellinger discontinuous at any $f_Q \in \mathcal{F}$.*

The discontinuity excludes the existence of estimators that have desirable properties in terms of unbiasedness and informativity [13].

An estimator $\hat{\theta}_n$ for $\theta$ is *locally unbiased* at $f_Q$ if there exists $\varepsilon > 0$ such that

$$\sup\{|E_G(\hat{\theta}_n) - \theta(f_G)| : f_G \in B(f_Q, \varepsilon)\} = 0,$$

where $E_G(\cdot)$ means taking expectation given $G \in \mathcal{Q}$ and where $B(f_Q, \varepsilon)$ is a ball,

$$B(f_Q, \varepsilon) = \left\{ f_G \in \mathcal{F} : \sum_{x=1}^{\infty} [f_Q^{1/2}(x) - f_G^{1/2}(x)]^2 \leq \varepsilon^2 \right\}.$$

The estimator $\hat{\theta}_n$ is *locally informative* at $f_Q$ if there exists $K(f_Q) > 0$ such that

$$\limsup_{\varepsilon \to 0} \sup\{E_G(\hat{\theta}_n^2) : f_G \in B(f_Q, \varepsilon)\} \leq K(f_Q).$$

An estimator (sequence) $\hat{\theta}_n$ for $\theta$ is *locally asymptotically unbiased* (*l.a.-unbiased*) at $f_Q$ with the rate of convergence $r(n) \geq n^{-1/2}$ if there exists $\varepsilon > 0$ such that

$$\lim_{m \to \infty} \limsup_{n \to \infty} \sup \left\{ \left| E_G \left[ l_m \left( \frac{\hat{\theta}_n - \theta(f_G)}{r(n)} \right) \right] \right| : f_G \in B(f_Q, \varepsilon n^{-1/2}) \right\} = 0,$$

where $l_m(z) = z - \operatorname{sign}(z) \max(|z| - m, 0)$. At $f_Q$, $\hat{\theta}_n$ is *locally asymptotically informative* (*l.a.-informative*) if there exist $\varepsilon > 0$ and $K(f_Q) > 0$ such that

$$\lim_{m \to \infty} \limsup_{n \to \infty} \sup \left\{ E_G \left[ l_m^2 \left( \frac{\hat{\theta}_n - \theta(f_G)}{r(n)} \right) \right] : f_G \in B(f_Q, \varepsilon n^{-1/2}) \right\} \leq K(f_Q).$$

THEOREM 3.1. *$\theta$ has no locally unbiased and locally informative estimator.*

THEOREM 3.2. *$\theta$ has no l.a.-unbiased and l.a.-informative estimator.*

Although bias is often the main concern, the discontinuity of $\theta$ will challenge our endeavor to reduce bias as a method for improving estimation accuracy.

THEOREM 3.3. *If $\{\hat{\theta}_{n,m}\}_{m=1}^{\infty}$ is a sequence of estimators for $\theta$ with fixed $n$, such that $\lim_{m \to \infty} |E_G(\hat{\theta}_{n,m}) - \theta(f_G)| = 0$ for $f_G$ in $B(f_Q, \varepsilon_0)$ with $\varepsilon_0 > 0$, then*

$$\lim_{m \to \infty} \sup \{ E_G(\hat{\theta}_{n,m}^2) : f_G \in B(f_Q, \varepsilon) \} = \infty, \qquad \varepsilon > 0.$$

Our ability to construct two-sided confidence intervals is also challenged. If, somewhere in $\mathcal{F}$, a confidence interval has a finite upper confidence limit with probability one, then, somewhere in $\mathcal{F}$, its coverage probability is zero [10].

THEOREM 3.4. *If $[\hat{\theta}_{n,l}, \hat{\theta}_{n,u}]$ is a confidence interval, then*

$$\sup \{ \Pr_Q(\infty \notin [\hat{\theta}_{n,l}, \hat{\theta}_{n,u}]) : Q \in \mathcal{Q} \} = 1$$

*implies that*

$$\inf \{ \Pr_Q(\theta(f_Q) \in [\hat{\theta}_{n,l}, \hat{\theta}_{n,u}]) : Q \in \mathcal{Q} \} = 0.$$

One can also consider the possibility that $\hat{\theta}_{n,u}$ is an upper confidence limit, that is,

(3.1) $\quad \inf \{ \Pr_Q(\hat{\theta}_{n,u} \geq \theta(f_Q)) : Q \in \mathcal{Q} \} \geq 1 - \alpha, \qquad \alpha \in (0,1).$

If the advertised confidence level is guaranteed, then the upper confidence limit will be infinite with large probability.

THEOREM 3.5. $\inf \{ \Pr_Q(\hat{\theta}_{n,u} = \infty) : Q \in \mathcal{Q} \} \geq 1 - \alpha.$

**4. Lower bounds.** We will construct lower bounds to $\theta$.

Although $\theta$ is discontinuous, it admits lower bounds, because of the following.

LEMMA 4.1. *$\theta$ is Kolmogorov lower semicontinuous at any $f_Q \in \mathcal{F}$.*

From [10], given $\varepsilon > 0$ and a distribution function $F_0$, with $F_Q(x) = \sum_{i=1}^{x} f_Q(i)$, the $\varepsilon$-lower envelope of $\theta$ at $F_0$ is

$$(4.1) \qquad \theta(F_0; \varepsilon) = \inf\{\theta(f_Q) : d(F_Q, F_0) \leq \varepsilon, f_Q \in \mathcal{F}\},$$

where $d(F_0, F_0^*)$ is the Kolmogorov distance of distribution functions $F_0$ and $F_0^*$,

$$d(F_0, F_0^*) = \sup\{|F_0(x) - F_0^*(x)| : x \in (-\infty, \infty)\}.$$

A conservative $1 - \alpha$ lower confidence limit is $\theta(\widehat{F}_n; \varepsilon_n)$, where $\widehat{F}_n(x) = \sum_{i=1}^{x} \hat{f}_n(i)$, $\hat{f}_n(x) = n_x/n$ and $\varepsilon_n$ is the $1 - \alpha$ quantile of the Kolmogorov distance of uniform$(0, 1)$ and the empirical distribution of $n$ random variables from it.

THEOREM 4.1. $\sup\{\Pr_Q(\theta(\widehat{F}_n; \varepsilon_n) \leq \theta(f_Q)) : Q \in \mathcal{Q}\} \geq 1 - \alpha$.

Calculating $\theta(\widehat{F}_n; \varepsilon_n)$ requires the solution of the optimization problem in (4.1). Given a grid $\{\xi_j\}_{j=1}^{J} \subset (0, \infty)$ with $Q = \sum_{j=1}^{J} \pi_j \delta(\xi_j)$, where $\delta(\lambda)$ is a degenerate distribution at $\lambda$, the discretized version of (4.1) is a linear program, due to the use of the Kolmogorov distance and the linearity of $\theta(f_Q)$ and $F_Q(x)$ in $Q$.

There are alternative lower bounds to $\theta$. Let $\mu(x) = \int \lambda^x \, d\Phi(\lambda)$ be the $x$th moment of a measure $\Phi$ over $(0, \infty)$ with $d\Phi(\lambda) = (e^\lambda - 1)^{-1} dQ(\lambda)$. Note that

$$\mu(0) = \theta(f_Q), \qquad \mu(x) = x! f_Q(x), \qquad x = 1, 2, \ldots.$$

When $\Gamma_k = (\mu(i+j))_{i,j=1}^{k}$ is positive definite, with $a_k = (\mu(j))_{j=1}^{k}$, define

$$(4.2) \qquad \theta_k = \theta_k(f_Q) = a_k' \Gamma_k^{-1} a_k.$$

THEOREM 4.2. *Let $\chi(Q)$ be the number of support points of $Q$. If $\chi(Q) < \infty$, then $\theta_1 < \cdots < \theta_{\chi(Q)} = \theta(f_Q)$, and $\theta_1 < \cdots < \lim_{k \to \infty} \theta_k = \theta(f_Q)$ otherwise.*

The approximation bias refers to $\theta_k - \theta$, whose absolute value decreases in $k$. The inferential challenge arises because the variance in $\theta_k$-estimation increases in $k$.

To find the condition under which the lower bound $\theta_k$ is Fisher consistent, we consider partitioning the mixture model $\mathcal{F}$ into "sieves" $\mathcal{F}_k = \{f_Q : \chi(Q) = k\}$.

THEOREM 4.3.   $\theta_k(f_Q) = \theta(f_Q)$ if $f_Q \in \mathcal{F}_k$; $\theta_k(f_Q) \leq \theta(f_Q)$ if $f_Q \in \bigcup_{j=k}^{\infty} \mathcal{F}_j$.

The lower bound $\theta_k$ is a functional that approximates $\theta$. A pre-existing nonparametric estimator for $c$ can also define an approximation functional to $\theta$. For example, from [6, 7, 9] one recognizes, with $s_i(f_Q) = \sum_{x=1}^{\infty} x^i f_Q(x)$,

$$\theta_{CB}(f_Q) = \frac{1 - f_Q(1)}{1 - f_Q(1)s_2(f_Q)/s_1^2(f_Q)} - 1,$$

$$\theta_{CL}(f_Q) = \frac{f_Q(1)\{s_2(f_Q) - s_1(f_Q)\} + s_1(f_Q)\{1 - f_Q(1)\}\{s_1(f_Q) - f_Q(1)\}}{\{s_1(f_Q) - f_Q(1)\}^2} - 1,$$

$$\theta_{DR}(f_Q) = \frac{1}{1 - f_Q(1)/s_1(f_Q)} - 1.$$

Unlike the $\theta_k$, it is not easy to specify the conditions under which each one is Fisher consistent, except that $\theta_{DR} = \theta_{CL} = \theta_{CB} = \theta$ when $Q = \delta(\lambda)$.

The lower bound $\theta_k$ is the odds of a mixing distribution from which the derived measure has the same first $2k + 1$ moments as $\Phi$ derived from $Q$.

THEOREM 4.4.   For $k \leq \chi(Q)$, there is a mixing distribution $Q_k$ with

$$\chi(Q_k) = k, \qquad \theta(f_{Q_k}) = \theta_k, \qquad f_{Q_k}(x) = f_Q(x), \qquad x = 1, 2, \ldots, 2k.$$

To estimate $\theta_k$, we consider the empirical moments $\hat{\mu}(x) = x!\hat{f}_n(x)$ and their matrices $\hat{a}_k = (\hat{\mu}(j))_{j=1}^k$ and $\widehat{\Gamma}_k = (\hat{\mu}(i+j))_{i,j=1}^k$. For $k \leq \hat{\chi}_n < \infty$, define $\hat{\theta}_k = \hat{a}_k' \widehat{\Gamma}_k^{-1} \hat{a}_k$, where $\hat{\chi}_n = \max\{k : |\widehat{\Gamma}_j| > 0, j = 1, 2, \ldots, k\}$.

THEOREM 4.5.   As $n$ goes to infinity, $\hat{\chi}_n$ estimates $\chi(Q)$ consistently when $\chi(Q) < \infty$. For finite $k \leq \chi(Q)$, as $n$ goes to infinity, $\hat{\theta}_k$ exists almost surely and $n^{1/2}(\hat{\theta}_k - \theta_k)$ converges to a zero-mean normal distribution.

Finally, an estimator for an approximation functional can also be calculated from $f_{\widehat{Q}}$ with $\widehat{Q}$ being the nonparametric MLE [12, 14]. Note that $\theta(f_{\widehat{Q}}) = \theta_{\chi(\widehat{Q})}(f_{\widehat{Q}})$ is the most greedy one among the $\theta_k(f_{\widehat{Q}})$ in terms of approximation bias reduction.

**5. Inference on $c$.**  We turn to unconditional inference on $c$.

As $c$ is identifiable given $\theta$ from $p_2(c, P)$ and $\theta$ is identifiable, it follows that $c$ is identifiable.

Let $\hat{c}_u$ be a $(1 - \alpha)$-level upper confidence limit for $c$, that is,

(5.1) $$\Pr_{c,P}(\hat{c}_u \geq c) \geq 1 - \alpha \qquad \forall c \geq 1, \forall P.$$

THEOREM 5.1. *For $(c, P)$, $\Pr_{c,P}(\hat{c}_u = \infty) \geq A(c, 1 - g_P(0)) - \alpha$, where*

$$A(c, \varrho) = \sum_{x=0}^{c} \min\left\{ \binom{c}{x} \varrho^x (1-\varrho)^{c-x}, \frac{e^{-c\varrho}(c\varrho)^x}{x!} \right\}.$$

The conclusion in Theorem 5.1 is slightly weaker than that in Theorem 3.5, as the distribution of $n$ retains a small amount information about $c$ from the testing affinity (see, e.g., [10]) of binomial$(c, \varrho)$ and binomial$(c', \varrho')$, with $\varrho = 1 - g_P(0)$, such that $c'\varrho'$ approaches $c\varrho$ when $c'$ goes to infinity. The bound $A(c, \varrho) - \alpha$ in Theorem 5.1 depends on $c$ and $P$ through the functional $\varrho$. From the fact that $A(c, 0) \equiv 1$, we can find a pair of $(c, \varrho)$ such that $A(c, \varrho)$ is arbitrarily close to one. For a fixed $c$, when the probability $\varrho$ of a class of being detected increases, the probability of the upper confidence limit being infinite will decrease. For an extremely large $\varrho$, one might have a negative value of $A(c, \varrho) - \alpha$. In particular, by Stirling's formula $c! \approx (2\pi c)^{1/2} (c/e)^c$, one has $A(c, 1) = e^{-c} c^c / (c!) \approx (2\pi c)^{-1/2}$. For $\alpha = 0.05$, $A(c, 1) > \alpha$ for $1 \leq c < 64$ and $A(c, 1) < \alpha$ for $c \geq 64$. Although the testing affinity $A(c, \varrho)$ is a function of both $c$ and $\varrho$, for $c$ relatively large it will change little in $c$ for a fixed $\varrho$. There exist lower bounds for $A(c, \varrho)$ that are functions of $\varrho$ only, for example, $A(c, \varrho) \geq 1 - 2^{-1} \varrho(1-\varrho)^{-1/2}$ [18].

Note that $\hat{c}_k = n(1 + \hat{\theta}_k)$ is a pseudo MLE for $c$ and is a consistent estimator for $c_k = c(1 + \theta_k)/(1 + \theta) \leq c$. In particular, $\hat{c}_1 = n + n_1^2/(2n_2)$ is given in [4]. The asymptotic variance of $\hat{c}_k$ increases in $c$, while that of $\log \hat{c}_k$ decreases in $c$ because both $c^{-1/2}(\hat{c}_k - c_k)$ and $c^{1/2}(\log \hat{c}_k - \log c_k)$ converge to zero-mean normal distributions as $c$ goes to infinity.

**6. Applications.** We consider two applications. The first (*cholera*) concerns an epidemic of cholera in a village in India [2, 15]. There were households affected by cholera but having no case. Note that $n_x$ is the number of households having $x$ cases, with $n_1 = 32$, $n_2 = 16$, $n_3 = 6$ and $n_4 = 1$ among $n = 55$ identified infected households with $S = 85$ cholera cases. The second (EST) concerns $S = 2586$ expressed sequence tags (ESTs) from which $n = 1825$ genes were found [14, 15]. An EST is a partial sequence identifying an mRNA and ESTs are generated by sequencing randomly selected clones in a cDNA library made from an mRNA pool. There were expressed genes from which no EST was generated. Note that $n_x$ is the number of expressed genes from which $x$ ESTs were generated, with $n_1 = 1434$, $n_2 = 253$, $n_3 = 71$, $n_4 = 33$, $n_5 = 11$, $n_6 = 6$, $n_8 = 3$, $n_x = 2$ for $x \in \{7, 10, 11, 16\}$ and $n_x = 1$ for $x \in \{9, 12, 13, 14, 23, 27\}$.

The estimates for approximation functions are shown in Table 1, together with the lower 5% quantiles of estimates from 400 model-based resamples, using the nonparametric MLE $\widehat{Q}$. All estimates are comparable in *cholera*, as $\chi(\widehat{Q}) = 1$. The pre-existing estimates are not comparable in EST, as $\chi(\widehat{Q}) > 1$. The linear program yields the conservative 95% nonparametric lower confidence limits: $\theta(\widehat{F}_n; \varepsilon_n) = 0.250$ with $n = 55$ and $\varepsilon_n = 0.180$ in *cholera*; $\theta(\widehat{F}_n; \varepsilon_n) = 1.408$ with $n = 1825$ and $\varepsilon_n = 0.032$ in EST. These bounds are considerably smaller than the resampling quantiles. If $\theta_1(f_{\widehat{Q}})$ is used to estimate $\theta$ in *cholera*, then a pseudo MLE for the number of infected households is 88. If $\theta_2(f_{\widehat{Q}})$ is used to estimate $\theta$ in EST, then a pseudo MLE for the number of expressed genes is 7392.

To learn something about the approximation bias, we treat $\widehat{Q}$ as the true distribution and read across the rows labeled $f_{\widehat{Q}}$ in Table 1, with the largest value of the $\theta_k$ being $\theta(f_Q)$. The other pre-existing approximation functionals are not better than the $\theta_k$ in EST because $\theta_{DR}/\theta = 0.41$, $\theta_{CL}/\theta = 1.46$ and $\theta_{CB} < 0$.

**7. Discussion.** Conditioning on the sample size $S$, the $Y_i$ arise from a multinomial distribution with index $c$ and probabilities $p_i = \lambda_i / \sum_{j=1}^{c} \lambda_j$. The multinomial model is more cumbersome analytically as the $Y_i$ are not independent. Just as in contingency table analysis using log-linear models, a Poisson-based analysis usually gives quantitatively similar or identical results, even for fixed size samples.

Results similar to those developed here can be established for a multiple-population species problem modeled by truncated mixtures of multivariate densities [16]. There are also lower bounds that can be developed for the total number of classes.

Estimating the population size by partially sampling a population is another important and difficult problem [5]. It could be investigated by means of various models of mixtures (e.g., binomial mixtures). Although the population size is nonidentifiable nonparametrically, we claim that by adapting and extending the techniques used here, we can show that confidence intervals for the population size must be one-sided, but identifiable lower bounds to the population size exist.

TABLE 1
*Estimates and the lower 5% empirical quantile of resample estimates $f_{\widehat{Q}}$*
*(cholera: 1st block, EST: 2nd block)*

|  | $\theta_{DR}$ | $\theta_{CL}$ | $\theta_{CB}$ | $\theta_1$ | $\theta_2$ | $\theta_3$ | $\theta_4$ | $\theta_5$ |
| --- | --- | --- | --- | --- | --- | --- | --- | --- |
| $\hat{f}_n$ | 0.593 | 0.544 | 0.484 | 0.582 | | | | |
| $f_{\widehat{Q}}$ | 0.608 | 0.608 | 0.608 | 0.608 | | | | |
| 5% quantile | 0.407 | 0.410 | 0.407 | 0.412 | | | | |
| $\hat{f}_n$ | 1.245 | 4.462 | −1.395 | 2.227 | 2.849 | 3.000 | 3.071 | 3.404 |
| $f_{\widehat{Q}}$ | 1.245 | 4.488 | −1.391 | 2.228 | 3.051 | 3.070 | 3.072 | 3.072 |
| 5% quantile | 1.120 | 3.222 | −1.755 | 1.964 | 2.432 | 2.446 | 2.455 | 2.455 |

## 8. Proofs.

PROOF OF LEMMA 2.1. Let $d\Psi(\lambda) = \lambda(e^\lambda - 1)^{-1} dQ(\lambda)$. As $\lambda e^{-\lambda} \leq e^{-1}$,

$$\{1 - g_P(0)\} \int e^{\lambda t} d\Psi(\lambda) = \int \lambda e^{-\lambda(1-t)} dP(\lambda) \leq (1-t)^{-1} e^{-1} \leq 1$$

for $t \leq 1 - e^{-1}$. The existence of a moment generating function (m.g.f.) implies that $\Psi$ is uniquely determined by its identifiable moments $\int \lambda^x d\Psi(\lambda) = (x+1)! f_Q(x+1)$, $x \geq 0$. The measure $\Psi$ and the distribution $Q$ are identifiable. □

The total variation distance $\tau(\psi, \phi)$ and the Hellinger distance $h(\psi, \phi)$ between two densities $\psi(x)$ and $\phi(x)$ over $\mathcal{R}^K$ with Borel field $\mathcal{B}$ are given by

$$\tau(\psi, \phi) = \int |\psi(x) - \phi(x)| = 2\sup\{|\Pr_\psi(B) - \Pr_\phi(B)| : B \in \mathcal{B}\},$$

(8.1)

$$h(\psi, \phi) = \left\{\int [\psi^{1/2}(x) - \phi^{1/2}(x)]^2\right\}^{1/2}.$$

Note that $\tau(\psi, \phi)$ and $h(\psi, \phi)$ satisfy

(8.2) $$h^2(\psi, \phi) \leq \tau(\psi, \phi) \leq 2h(\psi, \phi).$$

We introduce a useful single-parameter submodel of $\mathcal{F}$. Let $\pi(s)$ and $\eta(s)$ be two functions of $s \in (0, 1)$ with $\pi(s) \in (0, 1)$ and $\eta(s) \in (0, \infty)$. Given $Q$, define

(8.3) $$Q_s = (1 - \pi(s))Q + \pi(s)\delta(\eta(s)).$$

It is clear that

(8.4) $$\tau(f_{Q_s}, f_Q) = \sum_{x=1}^\infty |f_{Q_s}(x) - f_Q(x)| \leq 2\pi(s),$$

(8.5) $$\theta(f_{Q_s}) = (1 - \pi(s))\theta(f_Q) + \pi(s)\theta(f_{\eta(s)}).$$

LEMMA 8.1. Given $\varepsilon > 0$ and $f_Q \in \mathcal{F}$, $\omega(\varepsilon; \theta, f_Q, \mathcal{F}) = \infty$, where

$$\omega(\varepsilon; \theta, f_Q, \mathcal{F}) = \sup\{|\theta(f_Q) - \theta(f_G)| : f_G \in B(f_Q, \varepsilon)\}.$$

PROOF OF LEMMAS 3.1 AND 8.1. If $\pi^2(s) = \eta(s) = s^2$ in (8.3), then from (8.4) and (8.5), $\lim_{s \to 0} \theta(f_{Q_s}) = \infty$ and $\lim_{s \to 0} \tau(f_{Q_s}, f_Q) = 0$. By (8.2), one has $\lim_{s \to 0} h(f_{Q_s}, f_Q) = 0$ so that Lemmas 3.1 and 8.1 hold. □

PROOF OF THEOREMS 3.1 AND 3.3. Under Lemma 8.1, Theorem 3.1 and Theorem 3.3 hold because of Theorem 1 and Theorem 3 in [13], respectively. □

PROOF OF THEOREM 3.2. Let $\hat{\theta}_n$ be l.a.-unbiased and l.a.-informative for $\theta$ with the rate of convergence $r(n) \geq n^{-1/2}$. Let $s = 1/n$, $\pi(n^{-1}) = \varepsilon^2/(2n^2)$ and $\eta(n^{-1}) = 1/(r(n)n^3)$ in (8.3). Let $G_n = Q_{1/n}$ and

$$W_n = \frac{\hat{\theta}_n - \theta(f_Q)}{r(n)}, \qquad Z_n = \frac{\hat{\theta}_n - \theta(f_{G_n})}{r(n)}, \qquad d_n = \frac{\theta(f_{G_n}) - \theta(f_Q)}{r(n)}.$$

Note that $\lim_{n\to\infty} n\pi(n^{-1}) = 0$ and $\lim_{n\to\infty} d_n = \infty$, and that $f_{G_n} \in B(f_Q, \varepsilon n^{-1}) \subset B(f_Q, \varepsilon n^{-1/2})$ because $h^2(f_Q, f_{G_n}) \leq \tau(f_Q, f_{G_n}) \leq \varepsilon^2 n^{-2}$ from (8.2), (8.4) and (8.5). By investigating the proof of Theorem 2 in [13], with

$$u_{m,n} = 2\{E_Q[l_m^2(W_n)] + E_{G_n}[l_m^2(Z_n)] + 2E_Q|l_m(W_n)| \cdot d_n + d_n^2\},$$

due to the l.a.-informativeness and l.a.-unbiasedness, we have

(8.6) $|E_Q[l_m(Z_n)]|/d_n = 1 + o(1/d_n)$ as $n \to \infty$ and then $m \to \infty$,

(8.7) $|E_Q[l_m(Z_n)]|/d_n \leq |E_{G_n}[l_m(Z_n)]|/d_n + h(f_Q^{(n)}, f_{G_n}^{(n)}) \cdot u_{m,n}^{1/2}/d_n$.

Because $E_Q|l_m(W_n)| \leq E_Q^{1/2}[l_m^2(W_n)]$, by the l.a.-informativeness,

(8.8) $u_{m,n}/d_n^2 = 2 + o(1)$ as $n \to \infty$ and then $m \to \infty$.

For large $n$, from the proof of Lemma A.1 in [10], we have

$$\begin{aligned}(8.9) \quad h^2(f_Q^{(n)}, f_{G_n}^{(n)}) &= 2[1 - \{1 - h^2(f_Q, f_{G_n})/2\}^n] \\ &\approx nh^2(f_Q, f_{G_n}) \leq \varepsilon^2/n.\end{aligned}$$

By the l.a.-unbiasedness, from (8.7), (8.8) and (8.9), it follows that

$|E_Q[l_m(Z_n)]|/d_n = o(1)$ as $n \to \infty$ and then $m \to \infty$,

which is in contradiction to (8.6). This implies that Theorem 3.2 holds. □

PROOF OF THEOREM 3.4. Given $z > \theta(f_Q)$, let $\pi(s) = s$ and $\eta(s) = s/\{z - \theta(f_Q)\}$ in (8.3). As $\lim_{s\to 0} \tau(f_{Q_s}, f_Q) = 0$ and $\lim_{s\to 0} \theta(f_{Q_s}) = z$ from (8.4) and (8.5), $\{(f_Q, \theta(f_Q)): f_Q \in \mathcal{F}\}$ is dense in $\{(f_Q, z): f_Q \in \mathcal{F}, z \geq \theta(f_Q)\}$. Theorem 3.4 holds by applying Theorem 2.1 from [10]. □

PROOF OF THEOREM 3.5. Let $\pi(s) = s$ and $\eta(s) = s^2$ in (8.3). Because

$$\tau^2(f_{Q_s}^{(n)}, f_Q^{(n)})/8 \leq 1 - \{1 - \tau(f_{Q_s}, f_Q)/2\}^n$$
$$\leq 1 - (1 - s)^n$$

from Lemma A.1 in [10] and (8.4), we conclude that $\lim_{s\to 0}\tau(f_{Q_s}^{(n)}, f_Q^{(n)}) = 0$. From the condition in (3.1), the definitions in (8.1) and the fact that

$$|\Pr_Q(\hat{\theta}_{n,u} \geq \theta(f_{Q_s})) - \Pr_{Q_s}(\hat{\theta}_{n,u} \geq \theta(f_{Q_s}))|$$
$$\leq \sup\{|\Pr_Q(B) - \Pr_{Q_s}(B)| : B \in \mathcal{B}\},$$

we have by the triangle inequality,

$$\Pr_Q(\hat{\theta}_{n,u} \geq \theta(f_{Q_s})) + \tau(f_{Q_s}^{(n)}, f_Q^{(n)})/2 \geq \Pr_{Q_s}(\hat{\theta}_{n,u} \geq \theta(f_{Q_s})) \geq 1 - \alpha.$$

By letting $s$ go to zero, $\Pr_Q(\hat{\theta}_{n,u} = \infty) \geq 1 - \alpha$ as $\lim_{s\to 0}\theta(f_{Q_s}) = \infty$ from (8.5). This inequality holds for all $Q$, which implies that Theorem 3.5 holds. □

PROOF OF LEMMA 4.1. Let $Q$ and $G_m$ be in $\mathcal{Q}$ with $\lim_{m\to\infty}d(F_{G_m}, F_Q) = 0$. As a function of $f_Q(x)$, $x = 1, \ldots, 2k$, $\theta_k(f_Q)$ is continuous, so it is continuous in $F_Q$ on its domain. If $\theta_k(f_Q)$ exists, then $\theta_k(f_{G_m})$ will exist for sufficiently large $m$ and $\theta_k(f_Q) = \lim_{m\to\infty}\theta_k(f_{G_m}) \leq \liminf_{m\to\infty}\theta(f_{G_m})$. Because

$$\theta(f_Q) = \begin{cases} \theta_{\chi(Q)}(f_Q) \leq \liminf_{m\to\infty}\theta(f_{G_m}), & \chi(Q) < \infty, \\ \lim_{k\to\infty}\theta_k(f_Q) \leq \liminf_{m\to\infty}\theta(f_{G_m}), & \chi(Q) = \infty, \end{cases}$$

the odds $\theta(f_Q)$ is lower semicontinuous. □

PROOF OF THEOREM 4.1. This holds following application of (3.13) from [10] and Lemma 4.1. □

Write $M > 0$ if a matrix $M$ is positive definite. Given a sequence $(\mu(0), \mu(1), \ldots)$, define Hankel matrices $H_k = (\mu(i+j))_{i,j=0}^k$ and $\bar{H}_k = (\mu(i+j+1))_{i,j=0}^k$ for each $k$. The following summarizes some results in the Stieltjes moment problem.

LEMMA 8.2. *The sequence* $(\mu(0), \mu(1), \ldots)$ *of real numbers is the moment sequence of a measure* $\Phi$ *on* $(0, \infty)$ *if and only if* $|H_k| > 0$ *and* $|\bar{H}_k| > 0$ *for* $k < \chi(\Phi)$, *and, when* $\chi(\Phi) < \infty$, $H_k$ *and* $\bar{H}_k$ *have rank* $\chi(\Phi)$ *for* $k \geq \chi(\Phi)$.

PROOF OF THEOREM 4.2. Write $\Gamma_{k+1}$ and $\Gamma_{k+1}^{-1}$ as partitioned matrices,

$$\Gamma_{k+1} = \begin{bmatrix} \Gamma_k & b \\ b' & \mu(2k+2) \end{bmatrix}, \qquad \Gamma_{k+1}^{-1} = \begin{bmatrix} \begin{pmatrix} \Upsilon & v \\ v' & w \end{pmatrix} \end{bmatrix},$$

where $b = (\mu(k+2), \mu(k+3), \ldots, \mu(2k+1))'$, $w = (\mu(2k+2) - b'\Gamma_k^{-1}b)^{-1}$, $v = -w \cdot \Gamma_k^{-1}b$ and $\Upsilon = \Gamma_k^{-1} + w \cdot \Gamma_k^{-1}bb'\Gamma_k^{-1}$. Note that $|\bar{H}_k| = (-1)^k|\Gamma_k|(\mu(k+$

$1) - a_k' \Gamma_k^{-1} b$) because $\bar{H}_k$ can be obtained, by exchanging the rows $k$ times, from

$$\begin{bmatrix} \mu(k+1) & b' \\ a_k & \Gamma_k \end{bmatrix}.$$

As $|\Gamma_{k+1}| = |\Gamma_k|(\mu(2k+2) - b'\Gamma_k^{-1}b)$, it follows that $w = |\Gamma_k| \cdot |\Gamma_{k+1}|^{-1}$. Write

$$\begin{aligned}
\theta_{k+1} &= (a_k', \mu(k+1)) \cdot \Gamma_{k+1}^{-1} \cdot (a_k', \mu(k+1))' \\
&= a_k' \Upsilon a_k + 2\mu(k+1) a_k' v + w \cdot \mu^2(k+1) \\
&= a_k' (\Gamma_k^{-1} + w \cdot \Gamma_k^{-1} bb' \Gamma_k^{-1}) a_k \\
&\quad - 2w \cdot \mu(k+1) \cdot a_k' \Gamma_k^{-1} b + w \cdot \mu^2(k+1) \\
&= a_k' \Gamma_k^{-1} a_k + w \cdot (\mu(k+1) - a_k' \Gamma_k^{-1} b)^2 \\
&= \theta_k + |\Gamma_k| \cdot |\Gamma_{k+1}|^{-1} \cdot \{|\bar{H}_k|(-1)^{-k}|\Gamma_k|^{-1}\}^2.
\end{aligned}$$

This means that if $\theta_{k+1}$ exists, then $\theta_{k+1}$ and $\theta_k$ satisfy

$$\theta_{k+1} = \theta_k + |\bar{H}_k|^2 \cdot |\Gamma_k|^{-1} \cdot |\Gamma_{k+1}|^{-1}.$$

Note that $\bar{H}_k > 0$ when $\Gamma_{k+1} > 0$ so that $\theta_k < \theta_{k+1}$.

When $\chi(Q) < \infty$, write $|H_{\chi(Q)}| = |\Gamma_{\chi(Q)}| \cdot (\mu(0) - \theta_{\chi(Q)})$. From Lemma 8.2, $|H_{\chi(Q)}| = 0$, which means that $\mu(0) = \theta_{\chi(Q)}$ as $|\Gamma_{\chi(Q)}| > 0$. When $\chi(Q) = \infty$, write $|H_k| = |\Gamma_k|(\mu(0) - \theta_k)$. The sequence $\theta_k$ is strictly increasing in $k$ and bounded above by $\mu(0)$ so that $\xi = \lim_{k \to \infty} \theta_k$ exists. Consider $(\xi, \mu(1), \mu(2), \ldots)$ associated with Hankel matrices $H_{k,\xi}$ and $\bar{H}_{k,\xi}$. Note that $|\bar{H}_{k,\xi}| > 0$ because $\bar{H}_{k,\xi} = \bar{H}_k$, and $|H_{k,\xi}| > 0$ because $\theta_k < \xi$ and $|H_{k,\xi}| = |\Gamma_k|(\xi - \theta_k)$. From Lemma 8.2, $(\xi, \mu(1), \mu(2), \ldots)$ is a moment sequence of a measure $\Phi_\xi$ on $(0, \infty)$ with $\chi(\Phi_\xi) = \infty$. Let $d\Psi(\lambda) = \lambda \, d\Phi(\lambda)$ and $d\Psi_\xi(\lambda) = \lambda \, d\Phi_\xi(\lambda)$. Note that $\Psi$ and $\Psi_\xi$ have the same moment sequence and that $\Psi$ has an m.g.f. from the proof of Lemma 2.1. This implies that $\Psi = \Psi_\xi$, so that $\Phi = \Phi_\xi$ and $\xi = \mu(0)$. □

PROOF OF THEOREM 4.3. From Lemma 8.2, it follows that for $k \leq \chi(Q)$, $\Gamma_k > 0$ as $\Gamma_k$ is identical to the Hankel matrix $\bar{H}_{k-1}$ of the moments of a measure $\Psi$ with $d\Psi(\lambda) = \lambda \, d\Phi(\lambda)$. This observation and Theorem 4.2 imply that Theorem 4.3 holds. □

PROOF OF THEOREM 4.4. Let $H_{k,z}$ be obtained from $H_k$ with $\mu(0)$ replaced by $z \in \mathcal{R}$. If $\Gamma_k > 0$, then $|H_{k,z}| = |\Gamma_k|(z - \theta_k)$. When $\theta_k$ exists, because $\Gamma_i > 0$ and $\theta_i < \theta_k$ it follows that $|H_{i,\theta_k}| = |\Gamma_i|(\theta_k - \theta_i) > 0$ for $i = 1, \ldots, k - 1$. In addition, $|H_{k,\theta_k}| = 0$ and $\bar{H}_{k-1} > 0$. From [8], there exists a measure $\Phi_k$ with $\chi(\Phi_k) = k$ such that $\int d\Phi_k(\lambda) = \theta_k$ and $\int \lambda^x \, d\Phi_k(\lambda) = \mu(x)$, $x = 1, \ldots, 2k$. With $\Phi_k$ having no mass at zero, Theorem 4.4 holds by letting $Q_k = (e^\lambda - 1) \, d\Phi_k(\lambda)$. □

PROOF OF THEOREM 4.5. By the strong law of large numbers, the empirical moments, moment matrices and their determinants converge almost surely, implying the consistency of $\hat{\chi}_n$ and the almost sure existence of $\hat{\theta}_k$ for $k \leq \chi(Q)$ as $n$ goes to infinity. The delta method yields the asymptotic normality of $\hat{\theta}_k$ as $n^{1/2}(\hat{f}_{n,k} - f_{Q,k})$ converges to a multivariate normal distribution, where $f_{Q,k} = (f_Q(1), \ldots, f_Q(2k))'$ and $\hat{f}_{n,k} = (\hat{f}_n(1), \ldots, \hat{f}_n(2k))'$. □

PROOF OF THEOREM 5.1. With $\kappa$ given by (2.1) and $Q_s$ used to show Theorem 3.5, let $Q = \kappa(P)$, $P_s = \kappa^{-1}(Q_s)$ and $c_s$ be the integer part of $c\{1 + \theta(f_{Q_s})\}/\{1 + \theta(f_Q)\}$. Let $\tau_s = \tau(p_1(c_s, P_s), p_1(c, P))$. It can be shown that

$$c \cdot \frac{1 + \theta(f_{Q_s})}{1 + \theta(f_Q)} \in [c_s, c_s + 1),$$

$$\lim_{s \to 0} c_s = \infty,$$

$$\lim_{s \to 0} \frac{c_s}{1 + \theta(f_{Q_s})} = \frac{c}{1 + \theta(f_Q)},$$

$$|\tau_s - \tau(p_2(c_s, P_s), p_2(c, P))| \leq \sum_{n=0}^{c} \tau(p_3(n, P_s), p_3(n, P)) p_2(c, P).$$

Let $\varrho = 1 - g_P(0)$. Because the mean of $p_2(c_s, P_s)$ goes to that of $p_2(c, P)$ as $s$ goes to zero, $p_2(c_s, P_s)$ tends to a Poisson density with mean $c\varrho$,

$$\lim_{s \to 0} \tau(p_3(n, P_s), p_3(n, P)) = \lim_{s \to 0} \tau(f_{Q_s}^{(n)}, f_Q^{(n)}) = 0$$

and $\lim_{s \to 0} \tau_s = \lim_{s \to 0} \tau(p_2(c_s, P_s), p_2(c, P)) = 2 - 2A(c, \varrho)$. From the condition in (5.1) and the definitions in (8.1), and because

$$|\Pr_{c,P}(\hat{c}_u \geq c_s) - \Pr_{c_s,P_s}(\hat{c}_u \geq c_s)|$$
$$\leq \sup\{|\Pr_{c,P}(B) - \Pr_{c_s,P_s}(B)| : B \in \mathcal{B}\},$$

it follows by the triangle inequality that

$$\Pr_{c,P}(\hat{c}_u \geq c_s) + \tau_s/2 \geq \Pr_{c_s,P_s}(\hat{c}_u \geq c_s) \geq 1 - \alpha.$$

The proof is completed by letting $s$ go to zero. □

**Acknowledgments.** The authors thank the Editor, Associate Editor and referees for their useful comments which led to the improvement of this presentation.

Department of Statistics  
University of California  
Riverside, California 92521  
USA  
E-mail: cmao@stat.ucr.edu

Department of Statistics  
Pennsylvania State University  
University Park, Pennsylvania 16802  
USA  
E-mail: bgl@psu.edu